\documentclass[12pt]{amsart}
\usepackage{amsmath,amscd,amsthm,amsfonts, amssymb,amsxtra}
\usepackage[all]{xy}
\thanks{{\it 2000 MSC.} 55P42: Primary, 55P43. 
Secondary: 55P40, 55P65.}
\newtheorem{thm}{Theorem}[section]  
\newtheorem*{un-no-thm}{Theorem}
\newtheorem{cor}[thm]{Corollary}     
\newtheorem{lem}[thm]{Lemma}         
\newtheorem{prop}[thm]{Proposition}  
        
\newtheorem{bigthm}{Theorem}

\theoremstyle{definition}
\newtheorem{defn}[thm]{Definition}   

\theoremstyle{definition}
\theoremstyle{definition}
\newtheorem*{ques}{Question}
\theoremstyle{remark}
\newtheorem{rem}[thm]{Remark}

\newtheorem*{acks}{Acknowledgements}
 
\newtheorem*{out}{Outline}

\begin{document}
\title[Moduli of suspension spectra]{Moduli of suspension spectra}
\date{\today}
\author{John R. Klein}
\address{Wayne State University, Detroit, MI 48202}
\email{klein@math.wayne.edu}
\begin{abstract} For a $1$-connected spectrum $E$,
we study the moduli space of suspension
spectra which come equipped with a weak
equivalence to $E$. We construct
a spectral sequence converging
to the homotopy of the moduli space in positive degrees. 
In the metastable range,
we get a complete homotopical classification of
the path components of the moduli space.
Our main tool is Goodwillie's calculus of homotopy
functors.
\end{abstract}
\thanks{The author is partially supported by NSF Grant DMS-0201695}
\maketitle
\setlength{\parindent}{15pt}
\setlength{\parskip}{1pt plus 0pt minus 1pt}
\def\Top{\bold T\bold o \bold p}
\def\Sp{\bold S\bold p}
\def\vo{\varOmega}
\def\vs{\varSigma}
\def\smsh{\wedge}
\def\flush{\flushpar}
\def\id{\text{id}}
\def\dbslash{/\!\! /}
\def\codim{\text{\rm codim\,}}
\def\:{\colon}
\def\holim{\text{\rm holim\,}}
\def\hocolim{\text{\rm hocolim\,}}
\def\hodim{\text{\rm hodim\,}}
\def\hocodim{\text{hocodim\,}}
\def\Bbb{\mathbb}
\def\bold{\mathbf}
\def\Aut{\text{\rm Aut}}
\def\cal{\mathcal}
\def\frak{\mathfrak}

\section{Introduction}
A well-known property of 
the singular cohomology 
of a space is that it possesses
the structure of a commutative graded ring. 
There is no corresponding property for the singular
cohomology of
spectra. This discrepancy between spaces and spectra
is connected with the observation that a spectrum $E$ can fail to
have a  ``diagonal map'' ${E \to E\smsh E}$. 
Of course, a diagonal map exists when $E$ 
is a suspension spectrum. This motivates

\begin{ques} Let $E$ be a spectrum. When
can we find a based space $X$ and a weak equivalence of spectra
$\Sigma^{\infty} X \simeq E$? In how many ways?
\end{ques} 

The first result of this paper gives a criterion for deciding 
the existence part of the above question in the
metastable range.
Recall that  $E$ is
{\it $r$-connected} if its homotopy groups $\pi_*(E)$ vanish when $* \le r$.
Write $\dim E \le n$ if $E$ is, up to homotopy, a CW spectrum with
cells in dimensions $\le n$. Recall that the {\it $k$-th extended power}
$D_k(E)$ is the homotopy orbit spectrum of 
the symmetric group $\Sigma_k$ acting
on the $k$-fold smash product $E^{\smsh k}$.

\begin{bigthm}[Existence]\label{thmA} There is 
an obstruction
$$
\delta_E \in [E,\Sigma D_2(E)]\, ,
$$
which 
is trivial when $E$ has the homotopy type of a suspension spectrum.

Conversely, if $\delta_E = 0$,
$E$ is $r$-connected, $r \ge 1$ and $\dim E \le 3r{+}2$, then
$E$ has the homotopy type of a suspension spectrum.
\end{bigthm}

Before stating our second result, we comment on the
relation between Theorem \ref{thmA} and one of the early results of Kuhn
(\cite[Th.\ 1.2]{Kuhn-susp}) which says that a connected spectrum $E$
is a {\it retract} of a suspension spectrum if and
only if there is a weak equivalence
$$
\Sigma^\infty \Omega^\infty E \,\, \simeq \,\, \bigvee_{k \ge 1} D_k(E)\, . 
$$
Theorem \ref{thmA} implies one can remove the word
``retract'' from Kuhn's result in the metastable range. 
\medskip

Next consider the collection of pairs $(X,h)$ in which $X$ is 
a based space
and $h \: \Sigma^{\infty} X \to E$ is a weak (homotopy) equivalence.
Equate two such pairs $(X,h)$ and $(Y,g)$ if and only if there is a map
of spaces $f\:X \to Y$ such that $g\circ \Sigma^{\infty}f$ is homotopic
to $h$ (in particular, $f$ is
a homology isomorphism). This generates an equivalence relation.
Let $\Theta_E$ denote the associated set of equivalence classes.
\medskip

Our second result
identifies $\Theta_E$ in the metastable range.

\begin{bigthm}[Enumeration]\label{thmB}
Assume $\Theta_E$ is nonempty
and is equipped with basepoint. 
Then there is a basepoint preserving function
$$
\phi\:\Theta_E \to [E,D_2(E)]\,.
$$
If $E$ is $r$-connected, $r\ge 1$ and  $\dim E \le 3r+2$,
then $\phi$ is a surjection. If in addition $\dim E \le 3r+1$,
$\phi$ is a bijection.
In particular,
$\Theta_E$ has the structure of an abelian group when
$\dim E \le 3r+1$.
\end{bigthm}

Theorem \ref{thmB} leads to the possibility of 
calculating $\Theta_E$ in a number of simple cases.
For example, table 4.1 appearing in Mahowald's memoir \cite{Mahowald-memoir}
provides extensive 
calculations of $\Theta_E$ at the prime $2$ 
in the case of two cell spectra (see \S8 below).
\medskip

In fact, $\Theta_E$ is the
set of path components of a space 
$
\frak M_E
$  
which is defined as follows: let $C_E$ be the category whose objects are
pairs $(Y,h)$ with $Y$ a based space and $h\:\Sigma^{\infty}Y
\to E$ a weak equivalence. A morphism $(Y,h) \to (Z,g)$ consists
of a map of based spaces $f\:Y\to Z$ such that 
$g\circ\Sigma^{\infty}f = h$.
Then define
$$
\frak M_E \,\, := \,\, |C_E| \, ,
$$ 
i.e., the geometric realization of (the nerve of) $C_E$.\footnote{There are set-theoretic difficulties
presented by this definition, since $C_E$ isn't
a small category. The problem can be avoided in a variety of
ways, all furnishing the same homotopy type
(compare \cite[p.\ 379]{Wald_k_theory}).} We call 
$\frak M_E$ the {\it moduli space} of suspension structures
on $E$.
\medskip

Our third result gives a spectral sequence converging to the homotopy of
$\frak M_E$ in positive degrees. Formulating it requires some preparation.

For $q \ge 2$, let $W_q$ 
be the spectrum with $\Sigma_q$-action that classifies the
$q$-th layer of the Goodwillie tower of the identity functor
from based spaces to based spaces. Unequivariantly, $W_q$ is a wedge
of $(q{-}1)\!$ copies of the $(1{-}q)$-sphere spectrum (see 
Johnson \cite{Johnson-deriv}). If $q \le 1$, we take $W_q$ to be
the trivial spectrum. 

If $E$ is a spectrum, then
its $q$-fold smash power  $E^{\smsh q}$ has the structure of a spectrum
with (naive) $\Sigma_q$-action. Give the smash product $W_q \smsh E^{\smsh q}$
the diagonal $\Sigma_q$-action. We are then entitled to form
$$
W_q \smsh_{h\Sigma_q} E^{\smsh q} \, ,
$$
i.e, the homotopy
orbit spectrum of $\Sigma_q$ acting on $W_q \smsh E^{\smsh q}$. Let 
$F(E,W_q \smsh_{h\Sigma_q} E^{\smsh q})$ denote the function space
of spectrum maps from $E$ to $W_q \smsh_{h\Sigma_q} E^{\smsh q}$.
With these definitions we are ready to state our third result.
\medskip

\begin{bigthm}[Spectral Sequence] \label{thmC} Let $E$ 
be $1$-connected. Assume $\frak M_E$ is
non-empty and equipped with basepoint.
Then there is a first quadrant spectral sequence converging
to $\pi_{p{+}1}(\frak M_E)$ with 
$$
E^1_{p,q} := \pi_p(F(E,W_q \smsh_{h\Sigma_q} E^{\smsh q}))
$$
and  $d^1$-differential of bi-degree $(-1,1)$. 
\end{bigthm}
\medskip
Note when $\dim E \le  n$, we have
$E^1_{p,q} = 0$ for $p \le q + 1 - n$. Thus for
fixed $p$, $E^1_{p,q}$ is non-zero for finitely many $q$.
\medskip

Our last result says that the obstruction $\delta_E$
is $2$-local.

\begin{bigthm}[Localization] \label{localize}
For any connected spectrum $E$, the homotopy class
$\delta_E$  becomes trivial after inverting
the prime $2$.

Furthermore, if $f\:E \to E'$ is a map of spectra
which is
a $2$-local weak equivalence, then 
$\delta_E = 0$ if and only if $\delta_{E'} = 0$.

Consequently, if $E$ and $E'$ 
satisfy the connectivity assumptions of 
Theorem \ref{thmA}, 
then $E$ is weak 
equivalent to  a suspension spectrum
if and only $E'$ is weak equivalent to a suspension spectrum.
\end{bigthm}

\subsection*{Remarks} 
Since the paper of Berstein and Hilton \cite{B-H},
the problem of deciding when a space has the homotopy type of
a (single or iterated) suspension has been 
intensively studied. 

To the best of my knowledge, the literature 
contains  much less information
about the desuspension problem for spectra,
aside from the trivial stable range
case (Freudenthal's theorem) and
a $p$-local version considered by Gray
\cite{Gray-desuspend}.

My original interest in the desuspension question
for spectra came from embedding theory.
In a future paper, we intend to use the above results
to attack certain embedding questions.
\medskip

\begin{acks}
This paper would not have come into being had I
not had the fortune of discussing mathematics with Greg Arone,
Randy McCarthy and Tom Goodwillie in the middle to late
1990s. I am indebted to them for
their help. I also owe much to Bill Richter, who explained to
me the proof of the Berstein-Hilton theorem back in the 1980s.

I am grateful to Bob Bruner for
explaining to me how to compute
the stable homotopy of stunted projective space 
in low degrees, and for 
drawing my attention to the early work of Kuhn.
I would like to thank Nick Kuhn for some very helpful remarks,
especially concerning the proof of Lemma \ref{transfer}.
\end{acks}
\medskip

\begin{out} \S2 consists of preliminary material. 
The proof of Theorem \ref{thmA} is
contained in \S 3.  In \S4 we prove Theorem \ref{thmB}.
In \S5 we express the connected components of the moduli
space $\frak M_E$ as the classifying space
of a suitable  monoid. \S6 contains the proof of Theorem
\ref{thmC}. Theorem \ref{localize} is proved in \S7. In \S8 we provide examples
in connection with Theorem \ref{thmB}. In \S9 we discuss some loose ends.
\S10, which should be of independent
interest, outlines the classification of quadratic
homotopy functors from spectra to spectra.
The classification presented here
is due to Goodwillie, but I first learned about a 
special case of it from Kuhn.
The classification
gives a more concrete description of the
invariant $\delta_E$. 
\end{out} 
\medskip

\section{Preliminaries}

In this section we give the conventions and tools used 
throughout the paper. 

\subsection*{Spaces} All spaces will be compactly generated, and $\Top$
will denote the category of compactly generated spaces.
In particular, we make the convention that
products are to be retopologized with respect to the compactly
generated topology. 
Let $\Top_*$ denote the category of
based spaces. A {\it weak equivalence} of
spaces is shorthand for (a chain of) weak homotopy equivalence(s).

We use the usual connectivity terminology for spaces:
A non-empty space is {\it $r$-connected} if its homotopy
vanishes in degrees $\le r$ with respect to any choice of basepoint
(note: every non-empty space is $(-1)$-connected).
A map $A\to B$
of  spaces, with $B$ nonempty, is {\it $r$-connected} 
if for any choice of basepoint in $B$,
the  homotopy fiber with respect to this choice of basepoint
is an $(r{-}1)$-connected space. 

A commutative square of spaces
\begin{equation}
\begin{CD}
A @>>> C \\
@VVV @VVV \\
B @>>> D
\end{CD}
\label{cart-cocart}
\end{equation}
is said to be {\it $k$-cocartesian} if the map 
$\text{hocolim}(B \leftarrow A \to C) \to D$ is
$k$-connected. It is $\infty$-cartesian
if it is $k$-cartesian for all $k$. Similarly,
the square is {\it $k$-cartesian} if
the map $A \to \text{\holim}(B \to D \leftarrow C)$
is $k$-connected. The square
is $\infty$-(co)cartesian
if it is $k$-(co)cartesian for every $k$.

\subsection*{Spectra}
A {\it spectrum} $E$ will be taken to mean a 
collection of based spaces $\{E_i\}_{i \in {\Bbb N}}$
together with based (structure) maps $\Sigma E_i \to E_{i{+}1}$
where $\Sigma E_i$ is the reduced suspension of $E_i$.
A {\it morphism} of spectra $E \to E'$ consists of maps
$E_i \to E'_i$ that are compatible with the structure maps.
We denote the category of spectra by $\Sp$. 

A map of spectra is {\it $r$-connected} if it induces
a surjection on homotopy up through degree $r$ and
an isomorphism in degrees less than $r$. A spectrum 
is $r$-connected if the map to the zero object
(consisting of the one point space in each degree) is 
$(r{+}1)$-connected. 

A map of spectra is a {\it weak equivalence}
if it is $r$-connected for all integers $r$.
This notion of weak equivalence comes from a Quillen
model structure on $\Sp$.
In this model structure, the fibrant objects are
the $\Omega$-spectra (those spectra $E$ such that
$E_i \to \Omega E_{i{+}1}$ is a weak equivalence
for all $i$. A cofibrant object is 
(a retract of) a spectrum which is built up from the zero
object by attaching cells. The model
structure on $\Sp$ comes equipped
with functorial factorizations; in particular,
fibrant and cofibrant approximation is functorial.
For details, see Schwede \cite{Schwede}.

We typically apply fibrant and/or cofibrant approximation
functors to maintain homotopy invariance. To avoid clutter, we
usually suppress the application of these approximations
in the notation.

For example, if $E$ is a fibrant spectrum 
we write $\Omega^\infty E$
for the zero space $E_0$ of $E$. If $E$ isn't fibrant,
to get a good construction 
we first replace $E$ by its associated $\Omega$-spectrum $E^\sharp$
and then define $\Omega^\infty E$ as $E^\sharp_0$.

If $X$ is a based space, its {\it suspension spectrum} 
$\Sigma^\infty X$ has $j$-th space $\Sigma^j X$,
the $j$-fold reduced suspension of $X$
(for this to have the correct homotopy type, we assume
that $X$ is a cofibrant space; i.e., the
retract of a cell complex). In particular,
 $Q(X) := \Omega^\infty \Sigma^\infty X$ is
the reduced stable homotopy functor. A map
of spaces $Y \to \Omega^\infty E$ is adjoint to
a map of spectra $\Sigma^\infty Y \to E$.

The notion of $r$-cartesian and $r$-cocartesian
squares generalizes to spectra in the obvious way.
Note that an $r$-cartesian square is
the same thing as and $(r{+}1)$-cocartesian square
in the category of spectra. A commutative diagram
such as  \eqref{cart-cocart} but this time of spectra,
which is $\infty$-cartesian and with 
$C$ weak equivalent to the zero object, is called
a {\it fibration sequence}
(for spectra, a fibration sequence is the same thing
as a {\it cofibration sequence}). By slight abuse of notation,
we write the sequence as $A \to B \to D$.

As in the introduction,
we write $\dim E \le n$
if and only if $E$ is, up to homotopy, obtained from the zero object 
by attaching cells of dimension $\le n$. 
A spectrum is {\it finite} if it is built up from the zero object by
attaching a finite number of cells. 
A spectrum is {\it homotopy finite} if it is weak equivalent to
a finite spectrum. 

In this paper, we can get away with a notion of
smash product which is associative,
commutative and unital up to homotopy (see e.g.
Lewis, May and Steinberger \cite[Chap. 2]{LMS}).
We also need to know that  the
extended power spectrum 
$$
D_k(E) = (E^{\smsh k})_{h\Sigma_k} := 
(E\Sigma_k)_+\smsh_{\Sigma_k} E^{\smsh k} \, 
$$
exists, is functorial,  homotopy invariant and
coincides with the usual one in the case
of spaces. For details, see \cite[Chap.\ 6]{LMS}.

\subsection*{Truncation} 
Let $Y$ be a based space, $W$ be a $1$-connected
spectrum and $f\:\Sigma^\infty Y \to W$
a map of spectra. Assume $\dim W \le n$ and that $f$ is $n$-connected.

\begin{lem}\label{truncation} If $n \ge 0$, there exists a space $Z$, and an 
$(n{-}1)$-connected map $g\:Z \to Y$ such that the composite 
$$
\Sigma^\infty Z \overset{\Sigma^\infty g}\to \Sigma^\infty Y 
\overset f\to W
$$
is a weak equivalence.
\end{lem}

\begin{proof} If $n \le 1$, then $W$ is weak equivalent to the
zero object, and the proof is trivial in this case.

Assume then that $n \ge 2$.
Let $h\Bbb Z$ be the Eilenberg-Mac~Lane spectrum.
Then the induced map 
$$
h\Bbb Z \smsh \Sigma^\infty Y \to h\Bbb Z \smsh W
$$
can be thought of as a map of chain complexes (via the Dold-Kan
correspondence). Applying the {\it truncation lemma} \cite[4.1]{Klein},
We obtain a space $Z$ and  an $(n{-}1)$-connected map $Z \to Y$ such that
the composite
$$
h\Bbb Z \smsh \Sigma^\infty Z \to 
h\Bbb Z \smsh \Sigma^\infty Y \to h\Bbb Z \smsh W
$$
is a weak equivalence. But this composite is
obtained by smashing the composite
$$
\Sigma^\infty Z \to 
\Sigma^\infty Y \to W
$$
with $h\Bbb Z$. The result then follows by 
Whitehead's theorem.
\end{proof}

\begin{cor}[Freudenthal] 
\label{freud} Let $E$ be a an $r$-connected spectrum, $r \ge 1$.
Assume $\dim E \le 2r{+}2$. Then $E$ is weak equivalent to
a suspension spectrum.
\end{cor}

\begin{proof} This follows
from \ref{truncation} because
the map $\Sigma^\infty \Omega^\infty E \to E$
is $(2r {+}2)$-connected (see \ref{converge} below).
\end{proof}

\subsection*{Another description of the moduli space}
One can  use Lemma \ref{truncation} to give
an alternative description of the moduli space when $E$ is $1$-connected.
Let $C'_E$ be the (not full) subcategory of $C_E$ having the same objects, 
and whose morphisms
$(Y,h) \to  (Y',h')$ satisfy the condition that
$Y \to Y'$ is a {\it weak homotopy equivalence}. Let
$\frak M'_E$ denote the realization  of
$C'_E$.

\begin{prop} \label{m'=m} Assume that $E$ is
fibrant, cofibrant and $1$-connected. Then the inclusion
$
\frak M'_E \subset \frak M_E
$
is a homotopy equivalence.
\end{prop}

\begin{proof} Let $Y \mapsto Y^+$ denote 
the plus construction.
Applying functorial factorization, we can arrange it
so that the natural map $Y \to Y^+$ is
a cofibration.

If $(Y,h)$ is an object of $C_E$, 
then the $1$-connectedness of $E$ implies $Y$ has trivial first homology.
Therefore, $Y^+$ is a $1$-connected space. So if $(Y,h) \to (Z,h')$
is a morphism, the induced map of plus constructions
$Y^+ \to Z^+$ is a weak equivalence by Whitehead's theorem.

Since $E$ is fibrant and cofibrant, the space $\Omega^\infty E$
is cofibrant and $1$-connected. The natural map 
$\Omega^\infty E  \to (\Omega^\infty E)^+$ is a cofibration
which is a weak equivalence. Consequently, there is a retraction
$r\: (\Omega^\infty E)^+ \to \Omega^\infty E$.

Let $i\: C'_E \to C_E$ denote the inclusion.
Define a functor
$$
\pi\: C_E \to C'_E
$$ 
by
the rule $\pi(Y,h) = (Y^+,h^\flat)$, where 
$h^\flat\: \Sigma^\infty (Y^+) \to E$ is adjoint to the composite
$$
Y^+ \overset{(\hat h)^+} \to (\Omega^\infty E)^+
\overset r\to \Omega^\infty E \, . 
$$
Then there is an evident natural transformation
from the identity functor of $C_E$ to the composite
functor $i\circ \pi$. There is a similar evident natural
transformation from the identity functor of $C'_E$
to $\pi\circ i$. It follows that $i$ induces a homotopy
equivalence on realizations.
\end{proof}

\section{Proof of Theorem \ref{thmA}}

After defining the obstruction $\delta_E$,
the idea of the remainder of the proof will be to construct a 
highly connected map from a suspension spectrum
to $E$. The proof is then completed by applying Lemma
\ref{truncation}. 

\subsection*{Definition of the obstruction}
We describe below a certain fibration sequence of spectra
$$
D_2(E) \to S_2(E) \to E
$$
in which $S_2\: \Sp \to \Sp$ is a certain homotopy functor.
  
Assuming this construction has been specified, we have:

\begin{defn} The class $$\delta_E \in [E,\Sigma D_2(E)]$$ is the obstruction
to splitting the above fibration sequence, i.e., the homotopy
class of its connecting map to the right.
\end{defn}

\subsection*{Construction of the fibration} 
Consider the homotopy functor 
$$
\Sigma^\infty\Omega^\infty\:\Sp \to \Sp
$$
which assigns to a spectrum the suspension spectrum of its zero space.
Let $S_k(E)$ denote the $k$-th stage of the Goodwillie
tower of this functor and let $F_k(E) := \text{fiber}(S_k(E) \to S_{k{-}1}(E))$
denote the $k$-th layer.

The following result has been noted by several people, including 
Goodwillie, Arone, McCarthy, and 
Ahearn and Kuhn \cite[Cor.\ 1.3]{Ahearn-Kuhn}.

\begin{lem} \label{converge}
Assume $E$ is $r$-connected. Then the map
$$
\Sigma^\infty \Omega^\infty E \to S_k(E)
$$
is $((k{+}1)r + k{+}1)$-connected. Consequently, the
Goodwillie tower of $\Sigma^\infty \Omega^\infty E$
is convergent if $E$ is $0$-connected.
Furthermore, there is a natural weak equivalence of functors
$$
F_k(E) \,\,
 \simeq \,\,  D_k(E)\, .
$$
\end{lem}

Applying \ref{converge},  we see that the bottom of the tower
yields a fibration sequence of spectra
$$
D_2(E) \to S_2(E) \to E
$$
together with a $(3r{+}3)$-connected map $\Sigma^\infty\Omega^\infty E
\to S_2(E)$.

\begin{cor}   
Assume $E$ is $r$-connected and $\dim E \le 3r+3$.
Then $\delta_E = 0$ if and only if the map $\Sigma^\infty \Omega^\infty E
\to E$ admits a section up to homotopy.
\end{cor}

\begin{proof} The class $\delta_E$ is trivial if and only if 
the $S_2(E) \to E$
has a section. The dimension constraint
on $E$ and the lemma show that $S_2(E) \to E$
admits a section if and only if  $\Sigma^\infty \Omega^\infty E
\to E$ admits a section up to homotopy.
\end{proof}

\subsection*{} Now assume $E$ is $r$-connected, $r \ge 1$,  $\dim E \le 3r+3$
and $\delta_E = 0$. By the above remarks, we are entitled to choose a
homotopy section $\sigma \: E \to \Sigma^\infty \Omega^\infty E$.
Let $$\Omega^\infty\sigma\: \Omega^\infty E \to 
\Omega^\infty \Sigma^\infty \Omega^\infty E$$ be the 
corresponding map of zero spaces.

There is another map $c\: \Omega^\infty E \to 
\Omega^\infty \Sigma^\infty \Omega^\infty E$ ({\it not} an infinite
loop map) which is defined by taking the  adjoint
to the identity map of $\Sigma^\infty \Omega^\infty E$.

Let $Y$ be the homotopy pullback of the maps $\Omega^\infty \sigma$
and $c$. Thus we have an $\infty$-cartesian square
of spaces
\begin{equation}
\UseComputerModernTips 
\xymatrix{
Y \ar[rr]^{j} \ar[d]_{i}&& \Omega^\infty E \ar[d]^{\Omega^\infty \sigma}  \\
\Omega^\infty E \ar[rr]_{c} && \Omega^\infty \Sigma^\infty \Omega^\infty E 
}
\label{diag1}
\end{equation}
(commutative up to preferred homotopy).

\begin{lem} With respect to the
above assumptions, let $\hat j \: \Sigma^\infty Y \to E$ be the adjoint
to the map labeled $j$ in diagram (\ref{diag1}). Then $\hat j$
is $(3r+2)$-connected.
\end{lem}

\begin{proof} Consider the diagram of
spectra
\begin{equation}
\UseComputerModernTips 
\xymatrix{
\Sigma^\infty Y \ar[rr]^{\Sigma^\infty j} \ar[d]_{\Sigma^\infty i} 
&&\Sigma^\infty\Omega^\infty E \ar[rr]^{q_E} 
\ar[d]^{\Sigma^\infty \Omega^\infty \sigma} && E \ar[d]^{\sigma}\\
\Sigma^\infty \Omega^\infty E \ar[rr]_{\Sigma^\infty c\quad } 
&& \Sigma^\infty \Omega^\infty \Sigma^\infty \Omega^\infty E 
\ar[rr]_{\quad q_{\Sigma^\infty \Omega^\infty E}} 
&&\Sigma^\infty \Omega^\infty E 
}
\label{diag2}\end{equation}
in which the maps labeled with $q$ 
are the counits to the adjunctions for 
$E$ and $\Sigma^\infty \Omega^\infty E$.
The map $\hat j$ is therefore given by the
composite $q_E\circ \Sigma^\infty j$,
and the composite along the bottom,
$q_{\Sigma^\infty \Omega^\infty E} \circ \Sigma^\infty c$,
is clearly the identity. The left square commutes up to
a preferred homotopy. The right
square commutes on the nose.

The maps $\Omega^\infty \sigma$ and $c$ 
in diagram (\ref{diag1}) are both $(2r{+}1)$-connected.
Consequently, by the dual Blakers-Massey theorem
(see \cite{Good_calc2}) diagram (\ref{diag1})
is $(4r{+}3)$-cocartesian.
In particular, the left square in diagram (\ref{diag2}) is
a $(4r{+}3)$-cocartesian square of spectra.

Regarding the right square in diagram (\ref{diag2}), 
If $C$ denotes the homotopy cofiber of $\sigma$, then the evident map
$C \to D_2(E)$ is $(3r{+}3)$-connected. Similarly,
if $C'$ denotes the homotopy cofiber of $\Sigma^{\infty}\Omega^{\infty}\sigma$,
an argument using the Blakers-Massey theorem, which we omit, shows that
the evident map $C' \to \Sigma^\infty\Omega^\infty D_2(E)$ is
also $(3r{+}3)$-connected. 

Corollary \ref{freud} shows that the map 
$\Sigma^\infty\Omega^\infty D_2(E) \to D_2(E)$ is $(4r{+}4)$-connected
(because $D_2(E)$ is $(2r{+}1)$-connected).
We infer that the map $C' \to C$ is also $(3r{+}3)$-connected.
Therefore the right square in diagram (\ref{diag2}) is
$(3r{+}3)$-cocartesian.

Putting both squares together, it follows that 
diagram (\ref{diag2}) is $(3r{+}3)$-cocartesian.
Since the bottom composite is the identity map,
we infer that the top composite $\hat j$ is
$(3r{+}2)$-connected, as asserted.
\end{proof}

\subsection*{Completion of the proof} Assume
in addition to the above that
$\dim E \le 3r{+}2$. Since the map
$\hat j\: \Sigma^\infty Y \to E$ is $(3r{+}2)$-connected
we can apply \ref{truncation}.  This  gives a 
based space $Z$ and a based map $Z \to Y$ such that the  
composite $\Sigma^\infty Z \to \Sigma^\infty Y \to E$ is
a weak equivalence. The proof of 
Theorem \ref{thmA} is now complete.

\section{Proof of Theorem \ref{thmB}}

\subsection*{Step 1} Let $\cal S_E$ denote the set of homotopy 
classes of sections of the fibration $S_2(E) \to E$.
Note that the abelian group $[E,D_2(E)]$
acts freely and transitively on $\cal S_E$
(cf.\ Lemma \ref{section} below). Thus
if $\cal S_E$ is given a basepoint, it follows that 
$[E,D_2(E)]$ and $\cal S_E$ are isomorphic (the
isomorphism is dependent on the choice of basepoint).

There is a function
$$
\phi\:\Theta_E \to \cal S_E
$$
defined by sending $(Y,h)$ to the class represented by
$$
E \simeq \Sigma^\infty Y \to S_2(\Sigma^\infty Y) 
\simeq S_2(E)
$$
(the map in the middle is the preferred section
of  $S_2(\Sigma^\infty Y) \to \Sigma^\infty Y$).

If $E \to S_2(E)$ is a representative of $\cal S_E$,
and $\dim E \le 3r+2$, we can
perform the constructions in the previous section
to get an element of $\Theta_E$ (the element
is not necessarily unique). 
It is straightforward
to check that this element of $\Theta_E$ maps to the given
element of $\cal S_E$. Thus, when $\dim E \le 3r{+}2$,
the function $\phi\:\Theta_E \to \cal S_E$ is onto.

To complete the proof of Theorem \ref{thmB}, it
will be sufficient to show that $\phi$ is 
one-to-one when $\dim E \le 3r+1$. 
Choose a basepoint $(Y,h)$ for $\Theta_E$.
Then $\cal S_E$ inherits a basepoint and $\cal S_E$
becomes identified with $[E,D_2(E)]$. Thus we may
rewrite $\phi$ as a basepoint preserving function 
$$
\Theta_E \to [E,D_2(E)] \, ,
$$
where the basepoint of the codomain is the zero element.
Since $Y$ was chosen arbitrarily,
it suffices to show that
$\phi$ is one-to-one at the inverse image of the basepoint.
\medskip

\subsection*{Step 2}
We digress to develop a relative version of Theorem
\ref{thmA}. 
Suppose that $A \rightarrowtail E$
is a cofibration in the category of spectra.
We write $\dim(E,A) \le n$ if $E$ is obtained from $A$ 
up to homotopy by attaching
cells of dimension $\le n$. Assume that $A = \Sigma^\infty Z$
is a suspension spectrum. 
\medskip

\begin{thm}\label{relative} There is an obstruction
$$
\delta_{(E,A)} \in [E/A,\Sigma D_2(E)]
$$
whose triviality is necessary to finding a
cofibration $Z \to W$ and a weak equivalence 
$\Sigma^{\infty} W\simeq E$ extending the identity on $\Sigma^\infty Z = A$.
If $E$ is $r$-connected ($r \ge 1$) and $\dim(E,A) \le 3r{+}2$, 
then the vanishing of this obstruction is also sufficient.
\end{thm}

The obstruction $\delta_{(E,A)}$ is defined in the same
way as $\delta_E$ taking care to notice that the restriction of
$\delta_E$ to $A$ has a preferred trivialization.
The proof of \ref{relative} is virtually the same as the proof
of Theorem \ref{thmA}, so we omit it.

\subsection*{Step 3}
We now return to the proof of Theorem \ref{thmB}. 
We shall apply \ref{relative} in the following situation: 
choose a representative $(Y,h)$ for an element in $\Theta_E$
(call this the basepoint).
Suppose that $(Y',h')$ represents an another element
in $\Theta_E$. 
Both elements combine to give a weak equivalence 
$$
\Sigma^\infty (Y \vee Y') \overset\sim \to E\vee E \, .
$$
Applying \ref{relative}
to the pair $(E\smsh I_+,E\vee E)$ 
we get a necessary  obstruction
$$
\delta_{(E\smsh I_+,E\vee E)} \in [\Sigma E,\Sigma D_2(E)]
= [E,D_2(E)]
$$
 to finding
a cofibration $Y \vee Y' \rightarrowtail W$ of based spaces and a weak
equivalence 
$$
(\Sigma^\infty W,\Sigma^\infty (Y \vee Y')) \simeq (E\smsh I_+,E\vee E) \, .
$$  
It follows from this that the zig-zag $Y \to W \leftarrow Y'$ 
equates $(Y,h)$ with $(Y',h')$ in $\Theta_E$. 

Thus, fixing $(Y,h)$ and allowing $(Y',h')$ to vary defines
a function
$$
\Theta_E \to [E,D_2(E)]\, ,
$$ 
which is just another description of the function $\phi$
(we omit the details).

Assume that $E$ is $r$-connected and
$\dim E \le 3r{+}1$. Then \ref{relative}
shows that $(Y',h')$ maps to zero under $\phi$
if and only if $(Y',h')$ is the basepoint
of $\Theta_E$. This completes the proof
of Theorem \ref{thmB}.

\section{$\frak M_E$ as a classifying space}

Fix a cofibrant based space $Z$. 
Let $\Top_{*/Z}$ denote the category of based spaces
over $Z$. 
An {\it object} $y$ of this category consists of a based space
$Y$ together with a based map $p_Y\:Y \to Z$. 
A {\it morphism} $(Y,p_Y) \to (Y',p_{Y'})$
is given by a map of based spaces $f\:Y \to Y'$ 
such that $p_{Y'}\circ f = p_Y$.

Since $\Top_{*/Z}$ is an over category of the Quillen
model category $\Top_*$, it follows that
it too has the structure of a Quillen model category
(see \cite{Quillen}).
A {\it weak equivalence} is
a morphism $y \to y'$ whose underlying map of spaces is a
weak homotopy equivalence. We say $y \to y'$ is a {\it fibration}
if its underlying map of spaces is. We say that
$y \to y'$ is a {\it cofibration} if it satisfies the left lifting
property with respect to the acyclic fibrations.

Let $w\Top_{*/Z}$ denote the subcategory consisting of
the weak equivalences. Let $w_{(y)}\Top_{*/Z}$ denote the
full subcategory of $w\Top_{*/Z}$ consisting of
those objects connected to $y$ by a chain of
weak equivalences. 

The proof of the following proposition, which
we attribute to Waldhausen, is proved by the same method
as \cite[2.2.5]{Wald_k_theory}.
We omit the details.

\begin{prop}\label{wald-lemma} Assume that $y$ is
fibrant and cofibrant. Then
there is a weak equivalence
of spaces
$$
|w_{(y)}\Top_{*/Z}| \simeq BG(y) \, ,
$$
where $G(y)$ denotes the topological monoid of
self homotopy equivalences of $y$ in $\Top_{*/Z}$, and
$BG(y)$ is its classifying space.
\end{prop}
\medskip

We apply this in the following special case:
let $E$ be a $1$-connected fibrant and cofibrant spectrum.
The category $C'_E$, whose realization is the moduli space
$\frak M'_E$ (cf.\ \S2), is just the full subcategory of
$w\Top_{*/\Omega^\infty E}$ whose objects $Y$ are such that
the map $p_Y\:Y \to \Omega^\infty E$ is adjoint to a weak equivalence
of spectra. 

If we combine Proposition 
\ref{wald-lemma} with Proposition \ref{m'=m}, we 
obtain

\begin{cor} \label{m'=bg} Let $y$ be an object
of $C_E$ which is fibrant and cofibrant when considered as an object
of $\Top_{*/\Omega^\infty E}$.
Let $\frak M_{E,(y)}$ be the connected component of $\frak M_E$
which contains $y$. Then there is a homotopy equivalence
$$
\frak M_{E,(y)} \simeq BG(y) \, .
$$
\end{cor}

\section{Proof of Theorem \ref{thmC}}

\subsection*{Outline of the proof}

The basepoint $y = (Y,h)$ of $\frak M_E$ can be
taken as an object of $C_E$. 
Taking the plus construction if necessary,
we can assume without loss in
generality that $Y$ is $1$-connected (see the
argument in the proof of Proposition \ref{m'=m}).
We may also assume that $y$  is fibrant
and cofibrant when considered as
an object of $\Top_{*/\Omega^\infty E}$.

We will construct a tower of fibrations of based spaces
$$
\cdots \to T_3(y) \to T_2(y) \to  T_1(y) 
$$
such that 
\begin{itemize}
\item $T_1(y)$ is contractible;
\item for $k > 1$, there is a weak equivalence
$$
\text{fiber}(T_k(y) \to T_{k-1}(y))
\,\, \simeq \,\,  \Omega^\infty (W_k \smsh_{h\Sigma_k} Y^{\smsh k}) \, ,
$$
where  $W_k$ denotes the spectrum that classifies the
$k$-th layer of the Goodwillie tower of the identity
functor on based spaces.
\item there is a weak equivalence
$$
\Omega \frak M_{E,(y)} \simeq \lim_k T_y(y) \, .
$$
\end{itemize}
Assuming this has been done, we can define
the spectral sequence $\{E^r_{p,q}\}$ as the 
homotopy spectral sequence of  the tower $\{T_k(y)\}_k$.

We now digress to discuss generalities about section spaces
and the basic properties of the Goodwillie tower of the identity
functor.
\medskip

\subsection*{Digression} 
Suppose 
$p\:E \to Z$ is a fibration of based spaces.
We say that $p$ is {\it principal} if there
exists a commutative $\infty$-cartesian square
of based spaces
$$
\begin{CD}
E @>>> P \\
@Vp VV @VVV \\
Z @>>> B
\end{CD}
$$
such that $P$ is contractible.
\medskip

Suppose that
$Z$ is connected. If $p\: E\to Z$ is principal, there
is an ``action'' $\Omega B \times E \to E$. 
If there is a section $Z \to E$ one can combine it
with this action to produce a map
of fibrations $\Omega B \times Z \to E$
covering the identity map of $Z$. This implies that $p$ is
weak fiber homotopically trivial. Let $\text{Sec}(p)$
denote the {\it space of sections} of $p$. Then we have
shown

\begin{lem} \label{section} Assume $p\: E \to Z$ is principal.
Assume that $\text{\rm Sec}(p)$
is non-empty and comes equipped with basepoint.
Then there is a weak equivalence of based spaces
$$
\text{\rm Sec}(p)\,\, \simeq \,\, F(Z,\Omega B)\, .
$$
\end{lem}
\medskip

We next recall for the reader the basic properties
of the Goodwillie tower of the identity functor
on based spaces (cf.\ Goodwillie \cite{Good_calc1},
\cite{Good_calc2}, \cite{Good_calc3}, Johnson
\cite{Johnson-deriv}).

\begin{thm} 
\label{good-identity}
There is a tower of fibrations
of homotopy functors on based spaces
$$
\cdots \to P_2(X) \to P_1(X) 
$$
and compatible natural transformations
$X \to P_k(X)$ such that
\begin{itemize} 
\item $P_1(X) = Q(X)$ is the stable homotopy functor;
\item For $k\ge 2$, the fibration $P_k(X) \to P_{k{-}1}(X)$
is principal;
\item the $k$-the layer $L_k(X)$ is naturally
weak equivalent to the functor 
$$
X \mapsto \Omega^\infty (W_k \smsh_{h \Sigma_k} X^{\smsh k}) \,  ;
$$
\item if $X$ is $1$-connected, then the natural map
$$
X \to \lim_k P_k(X)
$$
is a weak equivalence.
\end{itemize}
\end{thm}

\subsection*{Completion of the proof} 
We are now in a position to define the tower $\{T_k(y)\}_k$.

\begin{defn} Let $T_k(y)$ be the space of lifts
$$
\UseComputerModernTips 
\xymatrix{
&  P_k(Y) \ar[d]\\
Y \ar[r] \ar@{-->}[ur] & P_1(Y)
}
$$
where $Y \to P_1(Y)$ is the natural map.
Note that $T_k(y)$ comes equipped with a basepoint defined by the
natural map $Y \to P_k(Y)$.
\end{defn}
\medskip

From the definition of $T_k(y)$, there is an evident tower of 
fibrations of based spaces
$$
\cdots \to T_2(y) \to T_1(y)
$$
with $T_1(y) = \ast$. Furthermore, the fiber
of the map $T_k(y) \to T_{k-1}(y)$ at the basepoint
is identified with the
space of lifts
$$
\UseComputerModernTips 
\xymatrix{
&  P_k(Y) \ar[d]\\
Y \ar[r] \ar@{-->}[ur] & P_{k{-}1}(Y)\, .
}
$$ 
Note that this is just
the space of sections of the pulled back fibration 
$Y \times_{P_{k{-}1}(Y)} P_k(Y) \to Y$. 
The latter fibration is 
principal and comes equipped with a preferred section.

Applying \ref{section} and \ref{good-identity}
we obtain
a weak equivalence of based spaces 
$$
\text{fiber}(T_k(y) \to T_{k{-}1}(y))\,\, \simeq  \,\,
F(Y,\Omega^\infty(W_k\smsh_{h\Sigma_k} Y^{\smsh k}))\, .
$$ 
Taking adjunctions, we see that
$F(Y,\Omega^\infty(W_k\smsh_{h\Sigma_k} Y^{\smsh k}))$
is weak equivalent to the space 
$F(E, W_k\smsh_{h\Sigma_k} E^{\smsh k})$ 
(where we are using the
identification $\Sigma^\infty Y \simeq E$).

To complete the proof of Theorem \ref{thmC}, it suffices
to identify the inverse limit of the tower $\{T_k(y)\}_k$.
It is clear that this inverse limit
is just the space
of lifts 
$$
\UseComputerModernTips 
\xymatrix{
&  \underset{k}\lim P_k(Y) \ar[d]\\
Y \ar[r] \ar@{-->}[ur] & P_1(Y)\, .
}
$$
Recall that $P_1(Y) = Q(Y)$ and
since $Y$ is $1$-connected, $\lim_k P_k(Y) \simeq Y$. With respect
to these identifications, we infer that this space of lifts
is weak equivalent to the underlying space of the topological
monoid $G(y)$. Applying \ref{m'=bg} and \ref{m'=m} completes the proof
of Theorem \ref{thmC}.

\section{Proof of Theorem \ref{localize}}

Let $\text{tr}\:D_2(E) \to E^{\smsh 2}$ be the transfer
and let $\pi \: E^{\smsh 2} \to D_2(E)$ be the projection
(\cite[Chap.\ 4]{LMS}).
It is well-known that the composite
$$
p\circ\text{tr}\: D_2(E) \to D_2(E)
$$
is a weak equivalence after inverting $2$.
To prove the first part of Theorem \ref{localize},
it is clearly enough to show

\begin{lem} \label{transfer} For a connected spectrum $E$, the class
$$
(\Sigma \text{\rm tr})\circ \delta_E \in [E,\Sigma E\smsh E]
$$
is trivial.
\end{lem}

\begin{proof} This is a consequence a result  
of Ahearn and Kuhn (see \cite[Cor.\ 1.7]{Ahearn-Kuhn}
in conjunction with the discussion after  
\cite[Cor.\ 1.2]{Ahearn-Kuhn}). For an alternative proof,
see the remark after
Corollary \ref{quadratic-classification} below.

In the case we need, 
Ahearn and Kuhn construct a 
natural transformation $\Phi\:S_2(E) \to E \smsh E$
fitting into a preferred 
homotopy commutative diagram  
$$
\UseComputerModernTips 
\xymatrix{
\Sigma^\infty \Omega^\infty E \ar[d]\ar[rr]^{\rm diagonal\qquad}
&& \Sigma^\infty (\Omega^\infty E  \smsh  \Omega^\infty E) \ar[d] \\
S_2(E) \ar[rr]_{\Phi} && E \smsh E \, .
}
$$
They then identify the composite map
$$
D_2(E) \to S_2(E) \overset{\Phi}\to E\smsh E
$$
with the transfer.\footnote{Randy McCarthy has informed me that 
that this result is also implicit in the
recent Ph.D.\ thesis of K.B. Bauer \cite{Baxter-Bauer}.}

Since $E \mapsto E^{\smsh 2}$ is a 
homogeneous quadratic functor, it coincides
with the second stage of its Goodwillie tower
and its first stage is trivial. So the map of 
Goodwillie towers
associated with $\Phi$
in degrees $\le 2$ may be displayed as
$$
\UseComputerModernTips 
\xymatrix{
D_2(E) \ar[dr] \ar[r]^{\rm tr}  & E\smsh E \ar@{=}[dr] \\
& S_2(E) \ar[r]\ar[d] & E\smsh E \ar[d] \\
& E \ar[r] &\,  \ast \, .
}
$$
If we take horizontal homotopy cofibers
of the maps from stage two to stage one, we obtain a
homotopy commutative diagram
$$
\UseComputerModernTips 
\xymatrix{
E \ar[r]\ar[d]_{\delta_E} & \ast \ar[d]\\
\Sigma D_2(E) \ar[r]_{\Sigma {\rm tr}} & \Sigma E \smsh E 
}
$$ 
which shows that $(\Sigma{\rm tr}) \circ\delta_E$ is
null homotopic. 
\end{proof}

We now prove the second part of Theorem \ref{localize}.
Let $f\:E \to E'$ be a map which is a $2$-local weak equivalence.
Assume that $E$ and $E'$ satisfy the connectivity hypotheses of
Theorem \ref{thmB}. Then the naturality 
of the obstruction shows that the diagram
$$
\begin{CD}
E @>\delta_E >>\Sigma D_2(E) \\
@Vf VV  @VV \Sigma D_2(f) V \\
E' @>>\delta_{E'} > \Sigma D_2(E')
\end{CD}
$$
is homotopy commutative. If we invert $2$, then
by the first part, the obstructions
$\delta_E$ and $\delta_E'$ vanish.
It we localize at 2, then $f$ and $\Sigma D_2(f)$ become
weak equivalences, and therefore $\delta_E$ and $\delta_{E'}$ coincide.
Thus, integrally, $\delta_E = 0$ if and only $\delta_E' = 0$.

\section{Illustrations of Theorem \ref{thmB}}

We give two examples.

\subsection*{Two cells} 
For $p > 1$ and $q \le 3p{-}3$, consider the suspension spectrum  
$E = \Sigma^\infty (S^p \cup_f e^{q{+}1})$ where $f\: S^q\to S^p$
is some map. According to Theorem \ref{thmB}, we have
$$
\Theta_E \,\, \cong \,\, [E,D_2(E)] \,\, \cong \,\, 
\pi_{q{+}1}(D_{2}(S^p)) \, ,
$$
where the second of these isomorphisms 
comes from that fact that a homotopy class $E \to D_2(E)$ 
maps the top cell of $E$ into $D_2(S^p) \subset D_2(E)$.
\medskip

The first non-trivial group occurs when $q = 2p{-}1$.
The group $\pi_{2p}(D_2(S^p))$ is isomorphic to
$\Bbb Z$ if $p$ is even and $\Bbb Z_2$ if $p$ is odd.
The distinct elements of $\Theta_E$ are represented by
the suspension spectra of $S^p \cup_g e^{2p}$,
with 
$$
g = f + k[\iota,\iota]\, ,
$$
 where
$k$ is an integer (mod 2 if $p$ is odd), and $[\iota,\iota]$
is the Whitehead product. 
The group structure on $\Theta_E$ is given by adding these
integers.
\medskip

In fact, at the prime 2, Mahowald has computed $\pi_{q{+}1}(D_2(S^p))$
for $q \le \min{(3p-3,2p + 29)}$ 
(see  \cite[table 4.1]{Mahowald-memoir},
see also Milgram \cite[table 13.5]{Milgram-unstable}). 
For example, assuming in addition $p \equiv 1 \text{ mod } 16$, the first few  
groups are
\medskip
$$
\vbox{
\offinterlineskip \tabskip=2pt
\halign{
\strut # &
\vrule # &
\hfil # \hfil &
\vrule # &
\hfil # \hfil &
\vrule # &
\hfil # \hfil &
\vrule # &
\hfil # \hfil &
\vrule # &
\hfil # \hfil &
\vrule # &
\hfil # \hfil &
\vrule # &
\hfil # \hfil &
\vrule # &
\hfil # \hfil &
\vrule #\cr
\omit& \multispan{17}{\hrulefill}\cr
& & $j$ & & $0$
& & 1& &  2 & & 3 && 4 && 5 && 6&\cr
\omit& \multispan{17}{\hrulefill}\cr
& & $\pi_{2p{+} j}(D_2(S^p))$ &&  $\Bbb Z_2$ & & $\Bbb Z_2$
& & $\Bbb Z_8$ & &  $\Bbb Z_2$& & $0$ && $\Bbb Z_2$
&& $\Bbb Z_{16} \oplus \Bbb Z_2 $
&\cr
\omit& \multispan{17}{\hrulefill}\cr
}}\, .
$$ 

\subsection*{Tori} Let $E = \Sigma^\infty (S^p \times S^p)$,
with $p > 1$.
Then 
$$
\Theta_E \,\, \cong \,\, [E,D_2(E)] \,\, \cong\,\,
 \pi_{2p}(D_2(S^p \vee S^p))\, .
$$
When $p$ is even, the group $\pi_{2p}(D_2(S^p \vee S^p))$ is isomorphic
to $\Bbb Z^{\oplus 3}$. When $p$ is odd, it is isomorphic
to $\Bbb Z_2 \oplus \Bbb Z_2 \oplus \Bbb Z$.
The elements of $\Theta_E$ can be represented by the suspension
spectra of the complexes
$$
(S^p \vee S^p) \cup_h e^{2p}
$$
with attaching map
$$
h = kx\circ [\iota ,\iota] + \ell y\circ [\iota,\iota] + m\omega \, ,
$$
where $\omega$ is the attaching map for the top cell of $S^p \times S^p$,
$x$ and $y$ are the summand inclusions $S^p \to S^p \vee S^p$,
and $k,\ell$ and $m$ are integers (take $k$ and $\ell$ modulo 2 if $p$ is odd).
Hence the elements of $\Theta_E$ are specified  by  triples $(k,\ell,m)$.
The identity element of $\Theta_E$ is $(0,0,1)$ and addition is 
given by
$$
(k,\ell,m) + (k',\ell',m') = (k+k',\ell +\ell', m+m'-1) \, .
$$

\section{Loose Ends}

\subsection*{Relation to the James-Hopf invariant}  
Another approach to desuspension questions
is to inductively desuspend cell-by-cell 
(with respect to a cell decomposition
of $E$). We explain here how this relates to our approach.

The idea is this:
when $E$ can be written
as a homotopy cofiber of a map of suspension spectra
$f\:\Sigma^\infty A \to \Sigma^\infty B$, it turns out that
$\delta_E$ is closely related to the James-Hopf invariant 
$$
H_2(f) \in [\Sigma^\infty A,D_2(\Sigma^\infty B)] \, .
$$
Explicitly, there is a homomorphism
$$
\psi\:[\Sigma^\infty A,D_2(\Sigma^\infty B)] \to 
[E,\Sigma D_2(E)]
$$
given by suspending, precomposing with the
connecting map $E \to \Sigma \Sigma^{\infty} A$ and postcomposing
with the inclusion $\Sigma D_2(\Sigma^\infty B) \to \Sigma D_2(E)$.
We assert that  $\psi$ maps $H_2(f)$ to $\delta_E$ (we defer
the proof to another paper). 

If we also assume that $B$ is $r$-connected ($r\ge 1$), $\dim B \le 2r{+}1$,
$A$ is $2r$-connected and $\dim A \le 3r{+}1$,
then we infer that $E$ is $r$-connected and $\dim E \le 3r{+}2$.
Furthermore, obstruction theory implies that $\psi$ is
an isomorphism. We conclude that $E$ desuspends if and only if 
$H_2(f) = 0$ (the `if' part of this statement is well-known).

\subsection*{Musings on the 
spectral sequence} It would be desirable to have a version of
the spectral sequence 
in Theorem \ref{thmC} which converges to $\pi_*(\frak M_E)$
in degree zero.
\medskip

The spectral
sequence still lacks a geometric interpretation.
Conjecturally, the spectral sequence should
be a packaging machine for the obstructions
to equipping a spectrum with the structure 
of an
\medskip

\centerline{``$E_\infty$-coalgebra
over the sphere spectrum.''}
\medskip

What I have in mind here 
comes from an (unproved) observation that $S_2(E)$ is a model for
$(E\smsh E)^{\Bbb Z_2} =$ the categorical fixed points of $\Bbb Z_2$ acting on 
$E \smsh E$,
and a choice of section
for the fibration
$
S_2(E) \to E
$
can be thought of as commutative (but not 
necessarily associative) ``diagonal'' for 
$E$.

It is tempting to conjecture that the spectral sequence
arises from a tower whose $k$-stage encodes, in some sense,
the moduli space of spectra equipped with choice of commutative diagonal
that is coherently homotopy associative up to order $k-1$.

\section{Appendix: On quadratic homotopy functors}

The invariant $\delta_E$
was defined as the connecting map associated with the
degree two homotopy functor $E \mapsto S_2(E)$. 
It is not difficult to classify the quadratic homotopy functors
from spectra to spectra. The classification, which is expressed
in terms of the Tate construction,
gives a more concrete way of understanding $\delta_E$.
I would like  to thank Tom Goodwillie and Nick Kuhn for explaining
these ideas to me. 
\medskip

Let $f\:\Sp \to \Sp$ be a homotopy functor of degree $\le 2$. 
This means that $f$ maps strongly cocartesian $3$-cubes to
$\infty$-cartesian (= $\infty$-cocartesian) $3$-cubes. We also
require $f$ to satisfy the usual technical conditions 
appearing in the functor calculus. 
Namely, we assume that $f$ preserves filtered homotopy colimits
(in particular, $f$ is determined up to weak equivalence
by its restriction to finite spectra) 
and is stably excisive (see \cite{Good_calc2} 
for the meaning of this language). We also assume
that $f$ is reduced in the sense that its value
at the zero object is contractible. 
\medskip

The general theory describes $f$ up to homotopy as sitting within a
fibration sequence
$$
W_2 \smsh_{h\Bbb Z_2} E^{\smsh 2} \to f(E) \to W_1 \smsh E\, ,  
$$
where $W_1$ is a spectrum and $W_2$ is a spectrum with $\Bbb Z_2$-action
(we remind the reader that $W_2 \smsh_{h\Bbb Z_2} E^{\smsh 2}$
denotes the homotopy orbits of $\Bbb Z_2$ acting diagonally
on $W_2 \smsh E^{\smsh 2}$).
The connecting map of this sequence is then a natural transformation
$$
u_f(E)\: W_1 \smsh E \to \Sigma (W_2 \smsh_{h\Bbb Z_2} E^{\smsh 2}) \, .
$$
Renaming $W :=W_1$ and $V := \Sigma W_2$. We can rewrite 
$u_f(E)$ as
$$
W \smsh E \to V \smsh_{h\Bbb Z_2} E^{\smsh 2} \, .
$$
Notice that $u_f(E)$ is 
a natural transformation from the linear
functor $E \mapsto W \smsh E$ to the homogeneous
quadratic functor $E \mapsto V \smsh_{h\Bbb Z_2} E^{\smsh 2}$,
and clearly, $f(E)$ is identified with the homotopy fiber of
$u_f(E)$. 
We will show how to classify such natural transformations.
\medskip

Let $\alpha$ denote the unique non-trivial 
one dimensional real representation of $\Bbb Z_2$.
If $j \ge 0$ is any integer, let $S^{j\alpha}$ be
the suspension spectrum of the one point compactification
of the direct sum of $j$-copies
of $\alpha$.  If $j < 0$, define $S^{j\alpha}$ as 
$F(S^{-j\alpha},S^0)$, the Spanier-Whitehead dual of 
$S^{-j\alpha}$. Note that $S^{j\alpha}$ has a $\Bbb Z_2$-action,
and there is an evident cofibration sequence of spectra with $\Bbb Z_2$-action
$$
S^{(j{-}1)\alpha} \to S^{j\alpha} \to (\Bbb Z_2)_+ \smsh S^j
$$
where the third term is given by inducing up the 
$j$-sphere with trivial action to a spectrum with free $\Bbb Z_2$-action.

Let us now consider $u_f$ applied to
$S^j\smsh E$. Then we have
$$
u_f(S^j\smsh E)\: W \smsh E\smsh  S^j \to 
(V\smsh E^{\smsh 2}) \smsh_{h\Bbb Z_2} (S^j \smsh S^{j\alpha})
$$
where we have taken the liberty of rewriting $S^j \smsh S^j$
with tranposition action as $S^j \smsh S^{j\alpha}$.
If we desuspend $j$-times, we get for each $j$ a map
$$
\Sigma^{-j}u_f(S^j\smsh E) \: W \smsh E  \to  
(V\smsh E^{\smsh 2})\smsh_{h\Bbb Z_2} S^{j\alpha}\, .
$$

\begin{lem} \label{u-compatible} These maps are compatible in the sense that
the diagram 
$$
\xymatrix{
W \smsh E\ar[rrr]^{\Sigma^{-j}u_f(S^j \smsh E)} 
\ar[drrr]_{\Sigma^{-(j{+}1)} u_f(S^{j{+}1}\smsh E)\quad }
&&&
(V\smsh E^{\smsh 2})\smsh_{h\Bbb Z_2} S^{j\alpha} \ar[d]\\
&&& (V\smsh E^{\smsh 2})\smsh_{h\Bbb Z_2} S^{(j{+}1)\alpha}  
}
$$
is homotopy commutative.
\end{lem}

\begin{proof} If we apply $u_f$ to the
evident cocartesian square of spectra
$$
\begin{CD}
S^j\smsh E @>i_+ >> D^{j{+}1} \smsh E\\
@V i_- VV @VVV \\
D^{j{+}1}\smsh E @>>> S^{j{+}1} \smsh E
\end{CD}
$$
we obtain a map of squares of spectra. This map induces a homotopy commutative
diagram
$$
\begin{CD}
\Sigma (W\smsh E \smsh S^j) @>\Sigma u_f(S^j\smsh E) >> 
\Sigma (V\smsh E^{\smsh 2} \smsh_{h\Bbb Z_2} (S^j \smsh S^{j\alpha})) \\
@V\simeq VV @VVV \\
W \smsh E\smsh  S^{j{+}1} @>> u_f(S^{j{+}1}\smsh E) >  
V \smsh E^{\smsh 2}\smsh_{h\Bbb Z_2} (S^{j{+}1} \smsh S^{(j{+}1)\alpha})
\end{CD}
$$
where the vertical maps are induced by maps from the pushouts
coming from $i_\pm$. The proof is now completed by desuspending
$j{+}1$ times.
\end{proof}

Recall that any spectrum $U$ with (naive) $\Bbb Z_2$-action comes
equipped with a {\it norm map} $U_{h\Bbb Z_2} \to U^{h\Bbb Z_2}$
from  homotopy orbits to homotopy fixed points (for details, see 
\cite{Adem-Cohen-Dwyer}, \cite[\S2]{WWII}, \cite{Greenlees-May},
\cite{klein;farrell-tate}). Recall that 
the {\it Tate spectrum} 
$$
U^{t\Bbb Z_2}  
$$
is the homotopy cofiber of the
norm map. Thus we have a {\it norm sequence}
$$
U_{h\Bbb Z_2} \to U^{h\Bbb Z_2} \to U^{t\Bbb Z_2}\, .
$$

\begin{lem} \label{tate} There is a weak equivalence
$$
U^{t\Bbb Z_2}  \simeq \Sigma \,\, (\underset{j}{\text{\rm holim }} 
U \smsh_{h\Bbb Z_2} S^{j\alpha})
$$
\end{lem}

\begin{proof} For each $j$ we have a cofibration sequence
$$
(U \smsh S^{(j{-}1)\alpha})^{t\Bbb Z_2} \to
(U \smsh S^{j\alpha})^{t\Bbb Z_2}
\to 
(U \smsh (\Bbb Z_2)_+ \smsh S^j)^{t\Bbb Z_2} \, .
$$
The spectrum $U\smsh (\Bbb Z_2)_+ \smsh S^j$
is the result of inducing $U \smsh S^j$ up to a 
free $\Bbb Z_2$-spectrum.
Consequently, its Tate spectrum is weakly contractible
(\cite[2.6]{WWII}). 
It follows that the map $(U\smsh S^{(j{-}1)\alpha})^{t\Bbb Z_2} \to
(U\smsh S^{j\alpha})^{t\Bbb Z_2}$ is a weak equivalence for all $j$.

Next, consider for each $j$ the norm sequence
$$
U\smsh_{h\Bbb Z_2} S^{j\alpha} 
\to 
(U \smsh S^{j\alpha})^{h\Bbb Z_2} 
\to 
(U \smsh S^{j\alpha})^{t\Bbb Z_2} \, .
$$
If we take the homotopy inverse limit over $j \in \Bbb Z$, the
middle term becomes weakly contractible (since homotopy
fixed points commute with homotopy inverse limits and 
the maps $U \smsh S^{j\alpha} \to  U \smsh S^{(j{+}1)\alpha}$
are unequivariantly null homotopic). The maps indexing
the homotopy limit
of the third term are all weak equivalences, so the
homotopy limit in this case is just 
$U^{t\Bbb Z_2} = (U \smsh S^0)^{t\Bbb Z_2}$. The result is
now immediate.
\end{proof}

Let $$u_f^\flat(E)\: W\smsh E \to 
\Sigma^{{-}1}(V\smsh E^{\smsh 2})^{t\Bbb Z_2}$$ be the 
natural transformation
 resulting from Lemmas \ref{u-compatible} and \ref{tate}.

\begin{cor} Up to homotopy, we have a factorization
$$
\xymatrix{ 
W\smsh E \ar[rr]^{u_f^\flat(E)\quad } \ar[drr]_{u_f(E)} 
&& \Sigma^{{-}1}(V\smsh E^{\smsh 2})^{t\Bbb Z_2}\ar[d]^\partial \\
&& 
V\smsh_{h\Bbb Z_2} E^{\smsh 2}
}
$$
where $\partial$ is the connecting map in the norm sequence.
\end{cor}

\begin{lem} When restricted to
finite spectra, the functor
$$
E \mapsto \Sigma^{{-}1} (V \smsh E^{\smsh 2})^{t\Bbb Z_2}
$$
has degree $\le 1$. 
\end{lem}

\begin{proof}(Sketch). There is a cofibration sequence
$$
\Sigma^{{-}1} (V \smsh E^{\smsh 2})^{t\Bbb Z_2} 
\to V \smsh_{h\Bbb Z_2}  E^{\smsh 2}
\to  (V \smsh  E^{\smsh 2})^{h\Bbb Z_2} \, .
$$
All functors in this
sequence have degree $\le 2$.
The functor in the middle is homogeneous of degree two,
and hence, its linearization vanishes. A direct calculation
which we omit shows that the homogeneous degree two part of
  $E \mapsto (V \smsh  E^{\smsh 2})^{h\Bbb Z_2}$ coincides
with  $E \mapsto V \smsh_{h\Bbb Z_2}  E^{\smsh 2}$
(this is where the assumption that $E$ ranges over finite
spectra is used, for I don't know whether the functor
$E \mapsto (V \smsh  E^{\smsh 2})^{h\Bbb Z_2}$
preserves filtered homotopy colimits).
We infer that 
$E \mapsto \Sigma^{{-}1} (V \smsh E^{\smsh 2})^{t\Bbb Z_2}$
has trivial homogeneous degree two part; it therefore has degree $\le 1$.
\end{proof}

A consequence of this lemma is that the natural transformation
$u_f^\flat$
is determined by its value at $S^0$: i.e.,
for finite spectra $E$ we have
$u_f^\flat(E) \simeq u_f^\flat(S^0) \smsh \id_E$. 
Since $u_f(E) \simeq \partial \circ u_f^\flat(E)$,
where $\partial$ is independent of $f$, it follows that
 the restriction of $u_f$ to finite spectra is
determined by $u^\flat_f(S^0)$. As
$f$ is determined
by its restriction to finite spectra,
we conclude

\begin{cor} \label{quadratic-classification}
The functor $f$  is determined up to weak
equivalence by the homotopy class
$$
[u_f^\flat(S^0)] \in [W ,\Sigma^{{-}1}V^{t\Bbb Z_2}] \, .
$$
\end{cor}
\medskip

\begin{rem} The above analysis
yields an alternative proof of Theorem
\ref{localize}. The point is that away from 2,
the connecting map $\partial$ in the norm sequence
is null homotopic.
\end{rem}
\medskip

We now apply \ref{quadratic-classification} in the case of the functor
$E \mapsto S_2(E)$. Here,
$W = S^0$ and $V = S^1$ (with trivial action). 
Furthermore, the natural transformation $u_f(E)$
is just $\delta_E$.

We see that $\delta_E$ is determined by a certain homotopy
class
$$
\delta^\flat \in
[S^0,(S^0)^{t\Bbb Z_2}]\, .
$$
With slightly more care, it is not difficult to
give an explicit construction: let
map $\iota\: S^0 \to (S^0)^{h\Bbb Z_2}$ be adjoint
to the evident unit map $(B\Bbb Z_2)_+ \to S^0$. Then $\delta^\flat$
represented by the composite
$$
S^0 \overset \iota \to  (S^0)^{h\Bbb Z_2} \to  (S^0)^{t\Bbb Z_2} \, .
$$
According to the Segal conjecture
\cite{C}, $\pi_0((S^0)^{t\Bbb Z_2})$
is isomorphic to the $2$-adic integers. The class $\delta^\flat$ 
represents a topological generator of this group.

\end{document}